\documentclass[10pt]{amsart}
\usepackage{amssymb,amsxtra}

\newcommand{\emp}{\emptyset}

\newcommand{\x}{\times}
\newcommand{\cs}{\mbox{$C^{*}$-algebra}}
\newcommand{\css}{\mbox{$C^{*}$-algebras}}

\newcommand{\N}{\mathbb{N}}

\newcommand{\Pos}{\mathbb{P}}

\newcommand{\ov}[1]{\mbox{$\overline{#1}$}}

\newcommand{\al}{\mbox{$\alpha$}}

\newcommand{\bt}{\mbox{$\beta$}}

\newcommand{\la}{\mbox{$\lambda$}}

\newcommand{\La}{\mbox{$\Lambda$}}

\newcommand{\si}{\mbox{$\sigma$}}

\newcommand{\Om}{\mbox{$\Omega$}}

\newcommand{\bc}{\begin{center}}

\newcommand{\ec}{\end{center}}
\newcommand{\be}{\begin{enumerate}}
\newcommand{\ee}{\end{enumerate}}
\newcommand{\beqn}{\begin{eqnarray}}
\newcommand{\eeqn}{\end{eqnarray}}
\newcommand{\beqns}{\begin{eqnarray*}}
\newcommand{\eeqns}{\end{eqnarray*}}
\newcommand{\bq}{\begin{quote}}
\newcommand{\eq}{\end{quote}}

\newcommand{\bi}{\begin{itemize}}
\newcommand{\ei}{\end{itemize}}
\newcommand{\bd}{\begin{description}}
\newcommand{\ed}{\end{description}}

\theoremstyle{plain}
\newtheorem{theorem}{Theorem}[section]
\newtheorem{lemma}[theorem]{Lemma}
\newtheorem{proposition}[theorem]{Proposition}
\newtheorem{corollary}[theorem]{Corollary}

\theoremstyle{definition}
\newtheorem{definition}[theorem]{Definition}

\theoremstyle{remark}

\numberwithin{equation}{section}

\begin{document}
\title[Tychonoff's theorem, path spaces]{Tychonoff's theorem for locally compact spaces and an elementary approach to the topology of path spaces}

\author{Alan L. T. Paterson}
\address{Department of Mathematics, University of Mississippi, University,
Mississippi 38677}
\email{mmap@olemiss.edu}
\author{Amy E. Welch}
\address{Department of Mathematics, University of Mississippi, University,
Mississippi 38677}
\email{amy3welch@yahoo.com}
\date{2003}
\begin{abstract}
The path spaces of a directed graph play an important role in the study of graph $\css$.  These are topological spaces that were originally constructed using groupoid and inverse semigroup techniques.  In this paper, we develop a simple, purely topological, approach to this construction, based on Tychonoff's theorem.  In fact, the approach is shown to work even for higher dimensional graphs satisfying the finitely aligned condition, and we construct the groupoid of the graph.  Motivated by these path space results, we prove a Tychonoff theorem for an infinite, countable product of locally compact spaces.  The main idea is to include certain finite products of the spaces along with the infinite product. We show that the topology is, in a reasonable sense, a pointwise topology.  
\end{abstract}

\maketitle

\section{Introduction}

Many of the most important $\css$ that have been studied in recent years (e.g. the Cuntz-Krieger algebras) arise in the following context.  We are given a directed graph $\mathcal{E}$
with sets of vertices and edges denoted by $V, E$ respectively.  Associated with the graph is a universal $\cs$ $C^{*}(E)$ and also a reduced $\cs$ $C^{*}_{red}(E)$.  There is now a large literature studying these graph $\css$ and the $\css$ of higher dimensional graphs.  A sample of papers in this area are: \cite{BPRS,CuntzSimple,CK,EL2,HR,Kumjian,KPask,KPR,KPRR,Pasks,Patbook,Patgraph,RS,RSYhigher,RSYthe,rg,Szy1,Szy2}.
One of the major approaches to graph $\css$, developed in a number of the papers above, is through the construction of a locally compact groupoid associated with the graph called the {\em graph groupoid}.  (Renault in his book \cite{rg} first developed the groupoid approach for the Cuntz graph (below) $\mathcal{E}_{n}$.) 

The graph groupoid is built up out of paths in $\mathcal{E}$ and indeed in all cases, the unit space of the groupoid is a space of paths.  In the row-finite case (i.e. when each vertex emits only finitely many edges) this space is the set $Z$ of all infinite paths in $\mathcal{E}$.  The space $Z$ is a zero-dimensional locally compact space.  While Tychonoff's theorem can, in the row-finite case, be easily applied to topologize $Z$, the case of non-row-finite graphs posed problems in the development of the theory.  Indeed, in Cuntz's fundamental paper (\cite{CuntzSimple}), the study of the $\cs$ $O_{\infty}$, corresponding to the graph $\mathcal{E}_{\infty}$ - an infinite bouquet of circles, not row-finite - requires a somewhat different treatment from that of the $\cs$ $O_{n}$ ($n$ finite) whose corresponding graph $\mathcal{E}_{n}$ is a bouquet of $n$ circles and is row-finite.  

Simple examples, such as $\mathcal{E}_{\infty}$, illustrate that in the non-row-finite case, $Z$ is not usually locally compact in the pointwise topology.  The problem of the existence of a locally compact groupoid associated with a non-row-finite graph was solved by Paterson in \cite{Patgraph} using the universal groupoid $G$ associated with the graph inverse semigroup $S$.  (See \cite[4.3]{Patbook} and \cite{KhoshSkan} for the construction of the universal groupoid of an inverse semigroup.)  A key step in this procedure is the construction of the unit space $X$ of $G$.  The important part of this unit space can be identified with 
$W=Y\cup Z$, where $Y$ is the set of finite paths of $\mathcal{E}$.  The space $W$ is then a locally compact metric space, and is used to construct the graph groupoid almost exactly as $Z$ was used in the row-finite construction.

It is natural to ask if it is possible to construct the topology on $W$ using purely topological means, with no recourse to inverse semigroups and their universal groupoids.  We will show that this can be done quite simply.  Indeed, more generally, this purely topological construction applies, not just to $1$-graphs, but even to $k$-graphs, and this is the setting of the second section.  We describe the groupoid for a general, finitely aligned, $k$-graph\footnote{In a lecture at the 2003 {\em Groupoidfest} meeting at the University of Colorado, Paul Muhly described joint work with C. M. Farthing and T. Yeend in which an inverse semigroup approach constructs the groupoid for finitely aligned, higher rank graphs, generalizing the $1$-graph case treated in \cite{Patgraph}.  The first author is grateful to Dr. Muhly for information about this work.} and in particular, the pointwise topology of the unit space $W$.  A related construction, in the context of partial group actions, has been given for Exel-Laca algebras (\cite{EL}).   

Motivated by the second section, we obtain, in the third section of the paper, an extension of the classical Tychonoff theorem to the case of locally compact, Hausdorff, second countable spaces.  It is elementary, of course, that a product of locally compact spaces $X_{i}$ is locally compact if and only if all but finitely many of the spaces are compact.  Our result (Corollary~\ref{cor:final}) requires no such restriction.  The crucial point is that we have to include, along with the infinite product $X_{1}\x X_{2}\x\ldots$, certain {\em finite} products of the spaces.  Precisely, we show that there is a natural locally compact metric topology on 
\[  X= (X_{1}\x X_{2}\x\ldots)\cup[\cup_{n=1}^{\infty}(X_{1}\x X_{2}\x\ldots X_{n})]. \]
The proof uses the one-point compactifications of the $X_{i}$'s and quotient topology results.

We have been unable to find this result in topology texts.  Something similar, though, appears in the theory of Stochastic processes, where adjunction of the one point at infinity (e.g. \cite[p.324]{Kallenberg}) gives an absorbing state, so that when a path enters into that state, it is ``dead'' and therefore can be regarded as a {\em finite} path.  However, we have been unable to find in that literature a detailed account of our result.  We also show that the topology on $X$ can be regarded as a ``pointwise'' topology as in the classical Tychonoff theorem.  It is, the authors think, intriguing  that the study of groupoids and $\css$, often regarded as {\em non-commutative} topology, suggests, in a natural way, a result in (commutative) topology itself!  

For simplicity, all topological spaces in this paper are locally compact, Hausdorff, second countable (and hence metrizable).  All graphs are directed and have countably many vertices and edges.  The first author is grateful to Marcelo Laca and Alex Kumjian for helpful discussions concerning the results of this paper.


\section{Path spaces and their topologies}

In this section, we will construct a locally compact metric topology on the path space of a finitely aligned $k$-graph $\La$.  This topology is that of pointwise convergence on the space $Y$ of finite paths in $\La$, and is remarkably simple to describe.  Local compactness follows from Tychonoff's theorem and we will describe the r-discrete groupoid $G_{\La}$ of $\La$ at the end of the section.  (The $\css$ of this groupoid are the $\css$ of the $k$-graph.)

We first discuss  higher rank graphs ($k$-graphs).  Such graphs $\La$ were introduced and studied by A. Kumjian and D. Pask in \cite{KPask}.  (Graphs of rank 1 are just directed graphs.)  They defined, in terms of generators and relations, the $\cs$ $C^{*}(\La)$ of $\La$.  Further, they showed, under the assumptions of row finiteness and the absence of sources, that the infinite path space $\La^{\infty}$ of $\La$  is a locally compact Hausdorff space.  Using shift maps $\si^{p}$, they obtained an r-discrete groupoid $\mathcal{G}_{\La}$ for $\La$, analogous to that for 1-graphs, and showed (among other results) that $C^{*}(\mathcal{G}_{\La})=C^{*}(\La)$.  We will show that their approach to the groupoid $\mathcal{G}_{\La}$ can be made to work when $\La$ is {\em finitely aligned}.  (So $\La$ need not be row finite and there may be some sources.)

We first recall some definitions (\cite{KPask}).  Let $\N=\{0,1,2,\ldots \}$ and $k\geq 1$.  Then $\N^{k}$ is ordered by: $m\leq n$ $\iff$ $m_{i}\leq n_{i}$ for $1\leq i\leq k$.  Also if $m,n\in \N^{k}$, then $m\vee n\in \N^{k}$ is given by: $(m\vee n)_{i}=\max\{m_{i}, n_{i}\}$.
Further, regard $\N^{k}$ as a small category with one object, $0$ as identity and addition as composition of morphisms.  A {\em $k$-graph} is a countable small category $\La$ together with a functor $d:\La\to \N^{k}$, the {\em degree map}, with unique factorization, i.e. if $\la\in \La$ and $m,n\in \N^{k}$ is such that $d(\la)=m+n$, then there are unique elements $\mu, \nu\in \La$ such that $\la=\mu\nu$ with $d(\mu)=m, d(\nu)=n$.  We write $\mu\leq \la$.  For $\la\in \La$, let $\La^{n}=d^{-1}(n)$; then 
$\La=\cup_{n\in \N^{k}} \La^{n}$.  Further $\La^{0}$ can be identified with the set of objects of $\La$.  Note that a $1$-graph $\La$ determines a directed graph with set of vertices $\La^{0}$ and set of edges $\La^{1}$, and $\La^{n}$ is just the set of paths of length $n$.  (The reader is referred to \cite{KPask,RS} for more details.)

In \cite{KPask}, paths in $\La$ are defined in terms of $k$-graph morphisms from the category
$\Om_{k}=\{(m,n)\in \N^{k}\x \N^{k}: m\leq n\}$ into $\La$.  For the purposes of the present paper, it is convenient to reformulate this in terms of sequences in $\La$ indexed by the elements of certain subsets of
$\N^{k}$.  (This reformulation is effectively contained in the account of \cite{KPask}.)  An advantage of this is that we can think of a path in $\La$ rather as we regard a path in an ordinary directed graph in terms of the sequence of its initial segments increasing with segment length. 

Let $N\in (\N\cup \infty)^{k}$.  Then an {\em $N$-path} in $\La$ is a sequence
$\{\la_{m}:m\in \N^{k}, m\leq N\}$ such that $d(\la_{m})=m$ and $\la_{m}\leq \la_{n}$ whenever $m\leq n$.  The range map on such a path is given by: $r(\{\la_{m}\})=r(\la_{0})$.
The $N$-paths are identifiable with the $k$-graph morphisms $x:\Om(N)\to \La$, where $\Om(N)$ is the small category whose set of objects is $O(N)= \{p\in \N^{k}: p\leq N\}$ and whose morphisms are pairs of the form $(m,n)\in \N^{k}\x \N^{k}$, where $m\leq n\leq N$, $r(m,n)=m, s(m,n)=n$ and $(m,n)(n,p)=(m,p)$.  Indeed, given an $N$-path $\{\la_{m}\}$, the morphism $x$ is obtained by factorizing $\la_{n}=\la_{m}x(m,n)$ where $d(\la_{m})=m, d(x(m,n))=n-m$, while conversely, given such a morphism $x$, then we just take $\la_{m}=x(0,m)$.  Let $W^{N}$ be the set of $N$-paths in $\La$ and 
\[  W=\cup_{N\in (\N\cup \infty)^{k}}W^{N}              \]
be the {\em path space} of $\La$.  If $N\in \N^{k}$, and $\{\la_{m}\}$ is an $N$-path, then by the uniqueness of factorization, we can identify the path with the single element $\la_{N}\in \La^{N}$.  So the set $Y$ of such paths can be called the {\em finite} path space, and can be identified with $\cup_{N\in \N^{k}}\La^{N}$.  The set 
\[   \cup_{N\in (\N\cup \infty)^{k}\setminus \N^{k}}W^{N}    \]
is called the {\em infinite} path space of $\La$, and is denoted by $Z$.  (Note that if
$N\in (\N\cup \infty)^{k}\setminus \N^{k}$, then $N_{i}=\infty$ for some $i$, and so 
$\{\la_{m}\}\in W^{N}$ has infinitely many components.)  So $W=Y\cup Z$.  Next, we need to be able to form products of paths $yw$ for $y\in Y, w\in W$ and $s(y)=r(w)$.
To define this let $y\in W^{N}, w\in W^{N'}$ and $n\in \N^{k}, n\leq N+N'$.  Choose $p\in \N^{k}$ such that 
$p\leq N', n\leq N+p$.  Then $(yw)_{n}$ is defined by: $yw_{p}=(yw)_{n}\mu$ where $d((yw)_{n})=n, 
d(\mu)=N+p-n$.  

Let $v\in \La^{0}$.  Define
\[  W_{v}=\{\{\la_{n}\}\in W: \la_{0}=v\}.   \]
Then $W=\cup_{v\in \La^{0}}W_{v}$.  We will show that, with a natural topology of ``pointwise convergence'', each $W_{v}$ is a compact Hausdorff space.  (So $W$ is a locally compact Hausdorff space with the disjoint union topology.)

We say that $\La$ is {\em finitely aligned} (\cite[Definition 5.3]{RS}) if given 
$\la, \mu\in Y=\La$, the set
\[  E_{\la,\mu}=\{\nu\in Y: \la\leq \nu, \mu\leq \nu, d(\nu)=d(\la)\vee d(\mu)\} \]
is finite.  As observed in \cite{RSYhigher}, finite alignment is the condition required for the new Cuntz-Krieger relation of \cite{RSYhigher} to be $C^{*}$-algebraic rather than spatial.  In the present paper, the condition has a very simple interpretation: it is exactly what is required to make the path space locally compact.  (For the $1$-graph case, finite alignment is trivially satisfied.)

We now define a map 
$\al:W\to 2^{Y}=\{0,1\}^{Y}$ by setting $\al(w)(y)=1$ if and only if $y\leq w$.  For each $a\in 2^{Y}$, let $A_{a}=\{y\in Y: a(y)=1\}$.  It is obvious that $\al$ is one-to-one, since if $\al(w)=a$, then for each $n$, $w_{n}$ is the unique element of $A_{a}$ in $W^{n}$. The next result identifies the image of $\al$ in $2^{Y}$.

\begin{lemma}   \label{lemma:inW}
Let $a\in 2^{Y}$.  Then $a\in \al(W)$ if and only if $A_{a}$ satisfies the following conditions:
\bi
\item[(i)] $A_{a}\ne \emp$;
\item[(ii)] if $y_{1}\leq y_{2}$ in $Y$ and $y_{2}\in A_{a}$, then $y_{1}\in A_{a}$;
\item[(iii)] $A_{a}$ is directed upwards, i.e. if $y_{1}, y_{2}\in A_{a}$ then there exists $y\in A_{a}$ such that $y_{1}\leq y, y_{2}\leq y$. 
\ei
\end{lemma}
\begin{proof}
Suppose that $a\in \al(W)$.  Then $a=\al(w)$ for some $w\in W$.  Then $A_{a}\ne \emp$ since 
it contains $w_{0}$.  (ii) follows since if $y_{1}\leq y_{2}$ and $y_{2}\leq w$, then $y_{1}\leq w$.  For (iii), suppose that $y_{1}, y_{2}\in A_{a}$.  Then $y_{1}\leq w, y_{2}\leq w$.  Let $y_{1}\in W^{N_{1}}, y_{2}\in W^{N_{2}}, w\in W^{N_{3}}$.  Then $N_{1}\leq N_{3}, N_{2}\leq N_{3}$.  (Note that the entries of $N_{1}, N_{2}$ are finite, but those of $N_{3}$ may not be.)
Let $N_{4}=N_{1}\vee N_{2}\leq N_{3}$ and $y=w_{N_{4}}$.  Then $y_{1}\leq y, y_{2}\leq y$ and $y\in A_{a}$.

Conversely, suppose that $a\in 2^{Y}$ satisfies (i), (ii) and (iii).  For each 
$y\in A_{a}$, let $N_{y}$ be such that $y\in W^{N_{y}}$.  Define $N$ to be the smallest element of 
$(\N\cup \infty)^{k}$ such that $N_{y}\leq N$ for all $y\in A_{a}$.  (So for each $i$, $1\leq i\leq k$, $N_{i}=\sup_{y\in A_{a}}(N_{y})_{i}$.) Let $n\in \N^{k}$, $n\leq N$.  Then for each 
$i$, there exists $y(i)\in A_{a}$ such that $n_{i}\leq (N_{y(i)})_{i}\leq N_{i}$.
By (iii), there exists $y\in A_{a}$ such that $y(i)\leq y$ for $1\leq i\leq k$.  Then
$n\leq N_{y}$.  Define $w_{n}=y_{n}\in \La$.  We claim that $w$ is well-defined, that 
$w\in W^{N}$ and that $\al(w)=a$.

To this end, suppose that $n\leq N'\in \N^{k}$, $N'\leq N$, and that $y'\in W^{N'}\cap A_{a}$.
Then by (iii), there exists $z\in A_{a}$ with $y\leq z, y'\leq z$.  Then 
$y_{n}=z_{n}=y'_{n}$ so that $w_{n}$ is well-defined. It also follows that $w\in W^{N}$.
We now show that $\al(w)=a$.  Suppose first that $a(y)=1$.  Then $N_{y}\leq N$ and by the above,
$y\leq w$ and $\al(w)(y)=1$.  Conversely, suppose that $\al(w)(y)=1$.  Then $y\leq w$.  By the above argument, there exists $y'\in A_{a}$ such that $N_{y'}\geq N_{y}$.  Then 
$y=w_{N_{y}}=y'_{N_{y}}\in A_{a}$ by (ii). Hence $\al(w)=a$.
\end{proof}  

\begin{theorem}                     \label{th:wvcompact}
Assume that $\La$ is finitely aligned.  Then for each $v\in \La^{0}$, the set $\al(W_{v})$ is a closed subset of $2^{Y}$.  So identifying $W_{v}$ with $\al(W_{v})$ under the map $\al$, $W_{v}$ is a compact metric space under the topology of pointwise convergence on $Y$. 
\end{theorem}
\begin{proof}
Let $\{w^{n}\}$ be a sequence in $W_{v}$ such that $\al(w^{n})\to a$ in $2^{Y}$.  It is sufficient to show that $A_{a}$ satisfies (i), (ii) and (iii) of Lemma~\ref{lemma:inW}.

(i) follows since $a(v)=\lim \al(w^{n})(v)=1$ since $v\leq w^{n}$ for all $n$.  For (ii), suppose that $y_{1}\leq y_{2}$ in $Y$ and $y_{2}\in A_{a}$.  Since $w^{n}\to a$ and $a(y_{2})=1$, it follows that $\al(w^{n})(y_{2})=1$ eventually.  By Lemma~\ref{lemma:inW}, $\al(w^{n})(y_{1})=1$ eventually, so that $a(y_{1})=1$, and (ii) follows.  Next, for (iii), suppose that $y_{1}, y_{2}\in A_{a}$.  We have to show that there exists $y\in A_{a}$ such that $y_{1}\leq y, y_{2}\leq y$.  To this end, again since $\al(w^{n})\to a$, we have eventually that $\al(w^{n})(y_{1})=1=\al(w^{n})(y_{2})$. We can assume that this is true for all $n$.  Let 
$N=d(y_{1})\vee d(y_{2})$ and let $N_{n}$ be such that $w^{n}\in W^{N_{n}}$.  Then 
$N_{n}\geq N$ for all $n$.  Let $y^{n}=w^{n}_{N}$.  Then 
$y^{n}\in E_{y_{1},y_{2}}$, and since $E_{y_{1},y_{2}}$ is, by finite alignment, finite, we can, by going to a subsequence, assume that $y^{n}=y$ is constant for all $n$.  Since 
$1=\al(w^{n})(y)\to 
a(y)$, we have $a(y)=1$ and $y\in A_{a}$.  So (iii) is satisfied and $a\in \al(W)$.  
\end{proof}
\begin{corollary}   \label{cor:wvcompact}
Let $\La$ be finitely aligned.  Then the path space $W$ is a locally compact metric space under the topology of pointwise convergence on $Y$.  This topology is the disjoint union topology of 
the compact spaces $W_{v}$.
\end{corollary}
\begin{corollary}   \label{cor:laW}
Let $\La$ be finitely aligned. Then the family of sets of the form $\la W_{s(\la)}$ forms a basis of compact open subsets for the topology of $W$.   
\end{corollary}
\begin{proof} 
This follows by intersecting $W$ by basic open subsets of $2^{Y}$.
\end{proof}

\noindent
{\bf Remarks for 1-graphs}\\

Let $\mathcal{E}$ be a directed graph and $E,V$ be the sets of edges and vertices respectively.  
Then $Y$ is the set of finite paths in $\mathcal{E}$ and $Z$ the set of infinite paths in $\mathcal{E}$.  Then $W=Y\cup Z$.  The graph $\mathcal{E}$ is said to be {\em row-finite} if for each $v\in V$, there are only finitely many edges $e$ starting at $v$.  For each $v\in V$, let $Z_{v}, Y_{v}$ be the sets of infinite, finite paths $z=z_{1}z_{2} \ldots $ of edges $z_{i}$ in $E$ starting at $v$.  Then $W_{v}=Y_{v}\cup Z_{v}$ and is a compact open subset of $W$.  If $\mathcal{E}$ is row-finite, then $Z_{v}$ is a closed subset of $W_{v}$, and for every $y\in W$, the singleton set $\{y\}$ is open in $W$.  

Let $\ell(w)$ be the length of the path $w$ ($0\leq \ell(w)\leq \infty$).  Abusing our earlier notation, let $w_{i}$ be the $i$th component of a path $w$. Then a sequence $w^{n}\to w$ in $W_{v}$ if and only if: either (a) $\ell(w)=\infty$ in which case $\ell(w^{n})\to \infty$ and $w^{n}_{i}\to w_{i}$ for all $i$, or (b) $\ell(w)<\infty$, in which case eventually, $\ell(w^{n})\geq \ell(w)$, $w^{n}_{i}\to w_{i}$ for $1\leq i\leq \ell(w)$ and $w^{n}_{\ell(w)+1}\to \infty$ (where $w^{n}_{\ell(w)+1}$ is taken to be $\infty$ if $\ell(w^{n})=\ell(w)$).  A similar result concerning convergent sequences can be formulated for finitely aligned $k$-graphs in general.  We will also see this kind of convergence again in the next section within a much more general, topological context.

A good illustration of the preceding comments (and indeed one of the main motivations for the theory developed in this paper) is provided by the Cuntz graph $\mathcal{E}_{\infty}$, which has one vertex $v$ and a countably infinite set of loops beginning and ending at $v$.  These loops will be labelled by the set $\Pos$ of positive integers $1,2,3,\ldots$.  Then 
\begin{equation} \label{eq:einfty}
W=\{v\}\cup(\cup_{k=1}^{\infty}\Pos^{k})\cup \prod_{n=1}^{\infty}\Pos =W_{v}
\end{equation}
is a compact metric space.  (For an alternative account of the topology on $W$ in this case, see \cite[pp.139-140]{rg}.)  Examples of convergent sequences in $W$ are: $123\ldots n\to 1234\ldots$, $12n111\ldots\to 12$, and  $nnn\ldots\to v$.

\vspace*{.1in}

Let $\La$ be a general, finitely aligned $k$-graph as above.  The path groupoid of $\La$ is then constructed in a manner analogous to that of \cite{KPask} but taking into account the presence of finite paths.   Let $G_{\La}$ be the set of triples of the form $(\la w,d(\la)-d(\mu),\mu w)$ for $\la, \mu\in \La$ and $w\in W$ where $r(w)=s(\la)=s(\mu)$.
Multiplication and inversion on $G_{\La}$ are defined in the usual way making it into a groupoid: 
$(x,n,y)(y,\ell,z)=(x,n+\ell,z)$ and $(x,n,y)^{-1}=(y,-n,x)$.  The family of sets
\[  Z(\la,\mu)=\{(\la w,d(\la)-d(\mu),\mu w): w\in W, r(w)=s(\la)=s(\mu)\}   \] 
is a basis for a topology on $G_{\La}$.  Indeed, if 
$(x,n,y)\in Z(\la,\mu)\cap Z(\la',\mu')$, then 
$(x,n,y)\in Z(\la\vee\la', \mu\vee\mu')\subset Z(\la,\mu)\cap Z(\la',\mu')$.
Further, $G_{\La}$ is a second countable, ample (\cite[p.48]{Patbook}), locally compact metric groupoid with counting measure as left Haar system (cf. \cite[Proposition 2.8]{KPask}).  For example, $Z(\la,\mu)$ is homeomorphic to the basic compact, open subsets 
(Corollary~\ref{cor:laW}) $\la W_{s(\la)}, \mu W_{s(\mu)}$ of the unit space $W$ of $G_{\La}$ under the range and source maps: $(\la w,d(\la)-d(\mu),\mu w)\to \la w$, 
$(\la w,d(\la)-d(\mu),\mu w)\to \mu w$.  It seems very likely (cf. \cite{KPask}) that the groupoid approach to the $\css$ of a $1$-graph (\cite{KPRR,Patgraph,rg}) can be extended to the case of finitely aligned $k$-graphs using the groupoid $G_{\La}$.


\section{Tychonoff's theorem for locally compact spaces}

The natural locally compact metric topology existing on $W_{v}$ for $\mathcal{E}_{\infty}$ ((\ref{eq:einfty})) motivates the search for a {\em locally compact} version of Tychonoff's theorem.  We will show in this section that the general version of this, for a product of locally compact spaces, holds. The topological techniques used to prove this generalization are quite different from those of the preceding section.  
 
So let $X_{n}$ be a second countable, locally compact metric space ($n=1,2,3,\ldots$).  Let $Y^{k}=X_{1}\x\ldots X_{k}$, $Y_{0}=\cup_{k=0}^{\infty}Y^{k}$, $Z=\prod_{n=1}^{\infty}X_{n}$ and 
\[                              W_{0}=Y_{0}\cup Z.                                   \]
In the preceding sentence, $Y^{0}$ is understood to be $\{0\}$.  ($0$ plays the role of the vertex $v$ in the $\mathcal{E}_{\infty}$ case.)  Let 
$Y=\cup_{k=1}^{\infty}Y^{k}$ and $W=Y\cup Z$.  So $W$ is obtained from $W_{0}$ by removing $0$.
(Our notation $W$ differs here from that of the preceding section (where vertices were contained in $W$) but we really want $W$ in the present context to involve only the $X_{i}$'s with no reference to $0$.)  It will be shown that there is a natural pointwise topology on $W_{0}$ which makes $W_{0}$ into a compact metric space.  Removing the point $0$ from $W_{0}$ then makes $W$ into a locally compact metric space. 

For each $n$, adjoin a point $\infty_{n}$ to $X_{n}$ to give the one-point compactification $X_{n}^{\infty}$ of $X_{n}$.  Then $X_{n}^{\infty}$ is a compact metric 
space (\cite[p.146, Exercise 4]{Nag}).
Let $A=\prod_{n=1}^{\infty}X_{n}^{\infty}$.  Then $A$ is a compact metric space by Tychonoff's theorem.  We now want to relate $A$ to $W_{0}$ in the hope that we can use the compactness of $A$ to obtain a ``pointwise'' compact topology for $W_{0}$.  Unfortunately, there are elements of $A$ that contain $\infty_{n}$'s as components, and yet the elements of $W_{0}$ do not contain any $\infty_{n}$'s.  The map $Q:A\to W_{0}$, that we now introduce, makes identifications in $A$ that solves this difficulty.

\begin{definition}
The map $Q:A\to W_{0}$ is defined as follows: 
\begin{equation*}          \label{eq:mapQ}
Q(\{x_{i}\})=
\begin{cases}
\{x_{i}\} \text{ if no $x_{i}$ is $\infty_{i}$}\\
(x_{1},x_{2},\ldots x_{n}) \text{ if there is a smallest integer $n\geq 1$ such that 
$x_{n+1}=\infty_{n+1}$}\\
0 \text{ if $x_{1}=\infty_{1}$}
\end{cases}
\end{equation*}
\end{definition}

It is obvious that $Q$ is onto $W_{0}$ and is 1-1 on $Z$.  We give $W_{0}$ the quotient topology.  For $x=\{x_{i}\}\in A$ let $N(x)$ be, for the 3 cases of the definition respectively, $\infty$, $n$ and $0$.  Then if $Q(x)=Q(y)$, it follows that $N(x)=N(y)$.  The map $Q$ is not 1-1 on $A$.

\begin{theorem}           \label{th:Wcompact}
$W_{0}$ is a compact metric space.
\end{theorem}
\begin{proof}
Let $R$ be the equivalence relation on $A$ determined by $Q$: so 
$\{x_{i}\}\sim \{x_{i}'\}$ if and only if $Q(\{x_{i}\})=Q(\{x_{i}'\})$.  Then $A/R=W_{0}$, and by \cite[p.105]{Bourbaki} and \cite[p.159]{Bourbaki2}, we have to show that $R$ (identified with its graph) is closed in $A\x A$.  To this end let $(x^{n},y^{n}) \to (x,y)\in A \x A$, where for all $n$, $(x^{n},y^{n})\in R$.  Write $x^{n}=\{x_{i}^{n}\}, y^{n}=\{y_{i}^{n}\}, x=\{x_{i}\}, y=\{y_{i}\}$. Then for all $i$, $\{x_{i}^{n}\}\to \{x_{i}\}, \{y_{i}^{n}\}\to \{y_{i}\}$.   Also, set $N(x^{n})=N(y^{n})=m^{n}$ for all $n$.  We have to show that $(x,y)\in R$.  There are two cases to be considered.

Suppose first that the sequence $\{m^{n}\}$ is bounded.  We can, by going to a subsequence, suppose that $m^{n}=m\geq 0$ is constant.   Then $x_{m+1}^{n}=\infty_{m+1}=y_{m+1}^{n}$ for all $n$, and so taking limits, $x_{m+1}=\infty_{m+1}=y_{m+1}$.  Next, for $1\leq j\leq m$, $x_{j}^{n}=y_{j}^{n}$ for all $n$ by the definition of $Q$, and taking limits, we get $x_{j}=y_{j}$ for $1\leq j\leq m$.  It follows that the first occurence of an $\infty_{k}$ in $x$ and $y$ is at exactly the same place, the $j$th place say, and $1\leq j\leq m+1$.  So 
$x_{i}=y_{i}$ for $1\leq i\leq j-1$, and $(x,y)\in R$.

Next suppose that the sequence $\{m^{n}\}$ is unbounded.  Taking a subsequence, we can, without loss of generality, suppose that $m^{n}\to \infty$.  Let $i\in \Pos$.  Since $m^{n}>i$ eventually, we have $x_{i}^{n}=y_{i}^{n}\in X_{i}$ for large enough $n$, and taking limits gives $x_{i}=y_{i}\in X_{i}^{\infty}$.  So $x=y$ and $(x,y)\in R$.
\end{proof}
\begin{corollary}    \label{cor:final}
The space $W=W_{0}\setminus\{0\}$ is a locally compact metric space in the relative quotient topology of $W_{0}$. 
\end{corollary}

For $w\in W_{0}$ let $\ell(w)$ be the length of $w$.  (We take $\ell(0)=0$.)  Note that for 
all $x\in A$, $N(x)=\ell(Q(x))$.
We now want to show that the topology on $W$ is, in a natural sense, a ``pointwise'' topology.  To this end, we need only check that the convergent sequences in $W$ are of the appropriate form.  The form for $W_{0}$ (and hence for $W$) is given in the following definition.

\begin{definition}  \label{def:conv}
Let $\{x^{n}\}$ ($x^{n}=\{x^{n}_{i}\}$) be a sequence in $W_{0}$ and $x\in W_{0}$.  We say that $x^{n}\to x$ {\em pointwise} if the following holds:
\bi
\item[(i)] if $x\in Z$ then $\ell(x^{n})\to \infty$ and $x^{n}_{i}\to x_{i}$ for all $i$;
\item[(ii)] if $x\in Y_{0}$, then $\ell(x^{n})\geq \ell(x)$ eventually, $x^{n}_{i}\to x_{i}$ for 
$1\leq i\leq \ell(x)$ and $x^{n}_{\ell(x)+1}\to \infty_{\ell(x)+1}$.
\ei
\end{definition}

(Earlier (\S 2), we met this kind of convergence for sequences in the path spaces of $1$-graphs.)
We now want to show that the convergent sequences for $W_{0}$ are exactly the pointwise convergent sequences, so that we can regard the topology of $W$ as a natural ``pointwise'' topology.  This will now be proved.  For notational ease, we write, for $x\in W_{0}$, $\ov{x}=Q^{-1}(\{x\})\subset A$.  This is an $R$-equivalence class in $A$.  

\begin{proposition}    \label{prop:converg}
$x^{n}\to x$ in $W_{0}$ $\iff$ every cluster point of every sequence $\{y^{n}\}$, where 
$y^{n}\in\ov{x^{n}}$, is in $\ov{x}$.
\end{proposition}
\begin{proof}
$\Rightarrow$ Suppose that for some sequence $\{y^{n}\}$, where $y^{n}\in\ov{x^{n}}$, there is a subsequence $y^{n_{k}}\to y\in A$.  Then $x^{n_{k}}=Q(y^{n_{k}})\to Q(y)=x$ so that $y\in \ov{x}$.\\
$\Leftarrow$ Suppose that $\{x^{n}\}$ does not converge to $x$.  By going to a subsequence, we can suppose that there is an open subset $U$ of $W_{0}$ with $x\in U$ and $x^{n}\in W_{0}\setminus U$ for all $n$.  Let $\{y^{n}\}$ be a sequence in $A$ with $y^{n}\in\ov{x^{n}}$.  
By compactness, there exists a cluster point $y$ of $\{y^{n}\}$.  Then for some subsequence $\{y^{n_{k}}\}$ of 
$\{y^{n}\}$,  we have $y^{n_{k}}\to y$.  Then $x^{n_{k}}=Q(y^{n_{k}})\to Q(y)$.  Since
$x^{n_{k}}$ is in the closed set $W_{0}\setminus U$, it follows that $Q(y)\in W_{0}\setminus U$ and $Q(y)\ne x$.  This contradicts $y\in \ov{x}$.
\end{proof}

\begin{theorem}   \label{th:final}
$x^{n}\to x$ in $W_{0}$ $\iff$ $x^{n}\to x$ pointwise.
\end{theorem}
\begin{proof}
There are two cases to be considered depending on whether $x\in Y_{0}$ or $x\in Z$.  We will give the proof for the first case, the proof for the second case being similar and easier.

Suppose then that $x\in Y_{0}$ and $x^{n}\to x$ in $W_{0}$.  Let $y^{n}\in \ov{x^{n}}$.  Then by 
Proposition~\ref{prop:converg}, every cluster point of $\{y^{n}\}$ is in $\ov{x}$. Suppose that it is not true that $\ell(x^{n})\geq \ell(x)$ eventually.  Then there is a subsequence $\{x^{n_{k}}\}$  such that $\ell(x^{n_{k}})=b<\ell(x)$.  So the $(b+1)$th entry of each $y^{n_{k}}$ is $\infty_{b+1}$.
We can suppose by compactness that $y^{n_{k}}\to y$ for some $y\in A$.  Since this convergence is pointwise in $A$, it follows that $y_{b+1}=\infty_{b+1}$ so that $y\notin \ov{x}$.  This is a contradiction.  So $\ell(x^{n})\geq \ell(x)$ eventually.  Now let $1\leq i\leq \ell(x)$.  We claim that $x^{n}_{i}\to x_{i}$.  For suppose not.  Then there exists a subsequence $\{x^{n_{k}}\}$ and a neighborhood $V_{i}$ of $x_{i}$ in $X_{i}$ such that for all $k$, $x^{n_{k}}_{i}\notin V_{i}$.
Let $y^{n_{k}}\in \ov{x^{n_{k}}}$.
By the above, without loss of generality, we can suppose that 
$N(y^{n_{k}})\geq \ell(x)$ for all $k$ and $y^{n_{k}}\to y$ for some $y$ in $A$.  Then 
$y_{i}=\lim x^{n_{k}}_{i} \in X^{\infty}\setminus V_{i}$.  But then 
$y\notin \ov{x}$ since the $i$th entry of $y$ does not equal that of $x$ and $i\leq \ell(x)$.  This is a contradiction.   So $x^{n}_{i}\to x_{i}$ for $1\leq i\leq \ell(x)$.  A similar argument shows that $x^{n}_{\ell(x)+1}\to \infty_{\ell(x)+1}$ so that by (ii) of Definition~\ref{def:conv}, $x^{n}\to x$ pointwise.

Conversely, suppose that $x^{n}\to x$ pointwise.  Let $y^{n}\in \ov{x^{n}}$ and $\{y^{n_{k}}\}$
be a subsequence such that $y^{n_{k}}\to y\in A$.  By Proposition~\ref{prop:converg}, we have to show that $y\in \ov{x}$.   To this end, we have that $y^{n_{k}}_{i}\to y_{i}$ in $X_{i}^{\infty}$ for all $i$.  We note that eventually, 
$N(y^{n_{k}})=\ell(x^{n_{k}})\geq \ell(x)$ and that
for all $1\leq i\leq \ell(x)$,  
$y^{n_{k}}_{i}=x^{n_{k}}_{i}\to x_{i}$ so that $y_{i}=x_{i}$ for $1\leq i\leq \ell(x)$.  It remains to show that $y_{\ell(x)+1}=\infty_{\ell(x)+1}$, for this would give that $y\in \ov{x}$. 

There are two cases to be considered.  Suppose first that $y^{n_{k}}_{\ell(x)+1}=\infty_{\ell(x)+1}$ for infinitely many $k$.  Then by pointwise convergence in $A$, $y_{\ell(x)+1}=\infty_{\ell(x)+1}$. For the second case, we can suppose, without loss of generality, that $y^{n_{k}}_{\ell(x)+1}\in X_{\ell(x)+1}$ for all $k$.  Then
$N(y^{n_{k}})=\ell(x^{n_{k}})\geq\ell(x)+1$.  So $x^{n_{k}}_{\ell(x)+1}=y^{n_{k}}_{\ell(x)+1}\in X_{\ell(x)+1}$, and 
$y_{\ell(x)+1}
=\lim_{k\to \infty} y^{n_{k}}_{\ell(x)+1}=\lim_{k\to \infty} x^{n_{k}}_{\ell(x)+1}=
\infty_{\ell(x)+1}$.  
\end{proof}

\noindent
{\bf Notes}  In the $\mathcal{E}^{\infty}$ case, i.e. where $X_{i}=\Pos$ for all $i$, there are two ways of defining a compact metric topology on $W_{0}$ ($=W_{v}$ ((\ref{eq:einfty}))).  The first of these, which we will call $(W_{0}, rel)$ is the relative topology obtained by regarding $W_{0}$ as a subset of $2^{Y_{0}}$ (using $\al$), discussed in the previous section.  The second is $(W_{0},quot)$ obtained by giving $W_{0}$ the quotient topology coming from the map $Q:A\to W_{0}$.  As one might expect, these are the same.  To see this, let $\bt:A\to 2^{Y_{0}}$ be given by: $\bt(x)(y)=1$ if $x$ begins with $y$ and is $0$ otherwise.  It is easy to check that $\bt$ is continuous and that 
$\bt=\al\circ Q$.  It follows by the definition of the quotient topology that the identity map from $(W_{0},quot)\to (W_{0},rel)$ is continuous.  An elementary compactness result then gives that the two topologies coincide on $W_{0}$.

In the case of Tychonoff's theorem for locally compact metric spaces (Theorem~\ref{th:Wcompact}), when the $X_{i}$'s are all the same space $X$, it is natural to ask if $W$ is the unit space of some locally compact groupoid.  (When $X=\Pos$, $W$ is the unit space of the Cuntz groupoid in the $\mathcal{E}_{\infty}$ case.)  There is a natural groupoid $\mathcal{G}_{X}$ whose elements are triples of the form $(yw,\ell(y)-\ell(y'),y'w)$ just as in the graph case.  However, we have been unable to topologize this groupoid in any reasonable way.


\end{document}